\documentclass[11pt]{article}
\usepackage{amsmath}
\usepackage{amssymb}
\usepackage{amsfonts}
\usepackage{mathrsfs}
\usepackage{shortvrb,psfrag}
\usepackage{enumerate}
\usepackage{color}
\usepackage[dvips]{graphicx}

\let\f=\frac

\let\Om=\Omega

\def\cD{{\cal D}}

\def\R{\Bbb R}

\def\no{\noindent}

\def\endproof{\hphantom{MM}\hfill\llap{$\square$}\goodbreak}

\newcommand{\beq}{\begin{equation}}
\newcommand{\eeq}{\end{equation}}
\newcommand{\ben}{\begin{eqnarray}}
\newcommand{\een}{\end{eqnarray}}
\newcommand{\beno}{\begin{eqnarray*}}
\newcommand{\eeno}{\end{eqnarray*}}

\newtheorem{Theorem}{Theorem}[section]
\newtheorem{Definition}[Theorem]{Definition}
\newtheorem{Proposition}[Theorem]{Proposition}
\newtheorem{Lemma}[Theorem]{Lemma}

\newtheorem{Remark}[Theorem]{Remark}

\begin{document}
\title{Non blow-up criterion for the 3-D Magneto-hydrodynamics equations in the limiting case}

\author{Wendong WANG\\[2mm]
{\small School of  Mathematical Sciences, Dalian University of Technology, Dalian 116024, P.R. China}\\[1mm]
{\small E-mail: wendong@dlut.edu.cn}}
\date{\today}
\maketitle

\begin{abstract}
In this paper, we prove that suitable weak solution $(u,b)$ of the 3-D MHD equations can be extended beyond $T$
if $u\in L^\infty(0,T;L^3(\R^3))$ and the horizontal components $b_h$ of the magnetic field satisfies the well-known Ladyzhenskaya-Prodi-Serrin condition.
\end{abstract}

\setcounter{equation}{0}
\section{Introduction}

We consider the 3-D incompressible Magneto-hydrodynamics (MHD) equations as follows:
\begin{align} \label{eq:MHD} \left\{
\begin{aligned}
&u_t-\Delta u+u\cdot \nabla u=-\nabla \pi+b\cdot \nabla
b,\\
&b_t-\Delta b+u\cdot \nabla b=b\cdot \nabla u ,\\
&\nabla\cdot u=\nabla\cdot b=0.
\end{aligned}
\right. \end{align}
Here $u$, $b$ describe the fluid velocity field and the magnetic field
respectively, $p$ is a scalar pressure.
The global existence of weak solution and local existence of strong solution to the MHD equations \eqref{eq:MHD}
were proved by Duvaut and Lions \cite{DL}. As the incompressible Navier-Stokes equations,
the regularity and uniqueness of weak solutions remains a challenging open problem. We refer to \cite{ST}
for some mathematical questions related to the MHD equations.

It is well-known that if the weak solution of the Navier-Stokes equations
satisfies the Ladyzhenskaya-Prodi-Serrin condition
\beno
u\in L^q(0,T;L^p(\R^3))\quad \textrm{with}\quad  \f 2 q+\f 3 p\le 1,\quad p\ge 3,
\eeno
then it is regular on $(0,T)\times \R^3$. Note that the limiting case $u\in L^\infty(0,T;L^3(\R^3))$
does not fall into the framework of energy method, which was proved by Escauriaza-Seregin-\v{S}ver\'{a}k \cite{ESS}.
 Wu \cite{Wu1,Wu2} extended Ladyzhenskaya-Prodi-Serrin type criterions to the MHD equations in terms of
 both the velocity field $u$ and the magnetic field $b$ for $p>3$.
However, some numerical experiments seem to indicate that the velocity field should play a more important role
than the magnetic field in the regularity theory of solutions to the MHD equations \cite{PP}.
Recently, He-Xin \cite{HX1} and Zhou \cite{Zhou1} have presented some regularity criterions to
the MHD equations in terms of the velocity field only.
Chen-Miao-Zhang \cite{CMZ1, CMZ2} extend and improve their results as follows:
if the weak solution of the MHD equations (\ref{eq:MHD}) satisfies
\beno
u\in L^q(0,T;B^s_{p,\infty})\quad \textrm{with}\quad\f 2 q+\f 3 p=1+s,\, \f 3 {1+s}<p\le \infty,\, -1<s\le 1,
\eeno
then it is regular on $(0,T)\times \R^3$.
Here $B^s_{p,\infty}$ is the Besov space. We refer to \cite{Cao, KL, WZ1} and references therein for
more relevant results. However, whether the condition on $b$ can be removed remains unknown in the limiting case (i.e., $(p,q,s)=(3,\infty,0$)).
The case $u,b\in L^\infty(0,T;L^3(\R^3))$ was considered by Mahalov-Nicolaenko-Shikin \cite{MNS}, and  Wang-Zhang \cite{WZ2} proved that
if $u\in L^\infty(-1,0;L^3(\R^3))$, and
\beno
 b\in L^{\infty}(-1,0;BMO^{-1}(\R^3)) \quad {\rm and}\quad b(t)\in VMO^{-1}(\R^3)\quad \textrm{for}\,\,t\in (-1,0].
\eeno
Then $(u,b)$ is H\"{o}lder continuous on $\R^3\times (-1,0]$.

Note that
the inclusion relation: $L^3(\R^3)\subsetneq VMO^{-1}(\R^3)$. It's interesting to ask: whether the condition on the magnetic field can be removed. 

Our main result is the following:
\begin{Theorem}\label{thm:uh serrin}
Let $(u,b)$ be a smooth solution of the MHD equations (\ref{eq:MHD})  in $(-1,0)\times \R^3$, which is also suitable. Assume that $u\in L^\infty(-1,0;L^3(\R^3))$ and $b$ satisfies one of the following conditions
\ben\label{eq:bh condition}
&&i)\nabla b_h\in L_t^qL_x^p((-\frac12,0)\times \R^3),\quad \frac3p+\frac2q=2,\quad \frac94\leq p<3;\nonumber\\
&&ii)\nabla b_h\in L_t^qL_x^p((-\frac12,0)\times \R^3),\quad \frac3p+\frac2q=2,\quad 3\leq p\leq \infty,\quad {\rm and}\\
&& \quad \quad b_h\in L_t^sL_x^l((-\frac12,0)\times \R^3),\quad \frac3l+\frac2s=1,\quad 9\leq l\leq \infty.\nonumber
\een
Then $(u,b)$ is regular on $\R^3\times (-1,0]$.
\end{Theorem}

\begin{Remark}
For simplicity, we assume that $(u,b)$ is smooth and suitable. At this moment, one can integrate by parts legitimately and the energy norms are finite, i.e.
\ben\label{eq:u b energy}
\|u\|_{L^\infty(-1,0;L^2(\R^3))\cap L^2(-1,0;\dot{H}^1(\R^3))}+\|b\|_{L^\infty(-1,0;L^2(\R^3))\cap L^2(-1,0;\dot{H}^1(\R^3))}<\infty.
\een
\end{Remark}

Motivated by the theory of Escauriza-Seregin-\v{S}ver\'{a}k \cite{ESS}, the main difficulty lies in proving that some necessary scaling invariant quantities of $u,b$ are bounded for compactness arguments. Using (\ref{eq:bh condition}), the vertical component $b_3$ of the magnetic field can be estimated by the energy method technically. On the other hand, the key is a class of careful interior regular criteria.
In fact, let $$G(b_h,p,q;r)\equiv r^{1-\frac3p-\frac2q}\|b_h\|_{L^q_tL_x^p(Q_r)},$$
where $Q_r=(-r^2,0)\times B_r$ and $B_r$ is a ball of radius $r$ centered at zero. We have the following more general interior criteria in the limiting case:
\begin{Theorem}\label{thm:uh interior}
Let $(u,b)$ be a suitable weak solution of the MHD equations (\ref{eq:MHD})  in $(-1,0)\times B_1$. Assume that $u\in L^\infty(-1,0;L^3(B_1))$, and $b$ satisfies the following conditions
\beno
&&i)\liminf_{r\rightarrow 0}G(b_h,p,q;r)=0,\quad\sup_{0<r<1}G(b_h,p,q;r)<\infty,\quad \frac3p+\frac2q=2,\quad 1\leq p\leq \infty;\\
&&ii)\sup_{0<r<1}G(b_3,l,s;r)<\infty,\quad \frac3l+\frac2s=2,\quad 1\leq s\leq \infty.
\eeno
Then $(u,b)$ is regular on $\R^3\times (-1,0]$.
\end{Theorem}

For $u\in L^{\infty}_tL_x^3$,  whether $b\in L^{\infty}_t(BMO^{-1}_x)$ or the above condition of $b$ without
$$\liminf_{r\rightarrow 0}G(b_h,p,q;r)=0$$ implies the regularity of $(u,b)$ is still unknown, where standard energy methods or backward uniqueness methods seem to be out of reach.



\setcounter{equation}{0}
\section{Preliminaries}

Let us first introduce the definitions of suitable weak solution.

\begin{Definition} Let $T>0$ and $\Om\subset \R^3$.
We say that $(u,b)$ is a suitable weak solution of the MHD equations (\ref{eq:MHD}) in $\Om_T=\Om\times (-T,0)$ if

(a)\,$(u,b)\in L^{\infty}(-T,0;L^2(\Om))\cap L^2(-T,0;H_0^1(\Om))$;

(b)\,$(u,b,\pi)$ satisfies the equations (\ref{eq:MHD}) in $\cD'(\Om_T)$;

%

(c)\,$\pi \in L^{\f 32}(\Om_T)$ and the following local energy
inequality holds: for a.e.\, $t\in [-T,0]$ \beno
&&\int_{\Om}(|u(x,t)|^2+|b(x,t)|^2)\phi dx+2\int_{-T}^t\int_{\Om}(|\nabla u|^2+|\nabla b|^2)\phi dxds\\
&&\leq
\int_{-T}^t\int_{\Om}\big[(|u|^2+|b|^2)(\triangle\phi+\partial_s\phi)+u\cdot\nabla\phi(|u|^2+|b|^2+2\pi)-(b\cdot
u)(b\cdot \nabla\phi)\big]dxds,\nonumber \eeno for any nonnegative
$\phi\in C_c^\infty(\R^3\times\R)$ vanishing in a neighborhood of
the parabolic boundary of $\Om_T$.
\end{Definition}


We define a solution $(u,b)$ to be regular at $z_0=(x_0,t_0)$ if $(u, b)\in L^\infty(Q_r(z_0))$ with
$Q_r(z_0)=(-r^2+t_0,t_0)\times B_r(x_0)$, and $B_r(x_0)$ is a ball of radius $r$ centered at $x_0$.
We also denote $Q_r$ by $Q_r(0)$ and $B_r$ by $B_r(0)$. For a function $u$ defined on $Q_r(z_0)$,
the mixed space-time norm $\|u\|_{L^{p,q}(Q_r(z_0))}$ is defined by
\beno
\|u\|_{L^{p,q}(Q_r(z_0))}^q:=\int_{t_0-r^2}^{t_0}\Big(\int_{B_r(x_0)}|u(x,t)|^pdx\Big)^{\f q p}dt.
\eeno

The following small energy regularity result is well-known, see \cite{HX1, MNS}. 

\begin{Proposition}\label{prop:small regularity}
Assume that $(u,b)$ is a suitable weak solution of (\ref{eq:MHD}) in $Q_1(z_0)$.
There exists an absolute constant $\varepsilon>0$ such that if
\beno
r^{-2}\int_{Q_r(z_0)}|u|^3+|b|^3+|\pi|^{\frac32}dxdt\leq \varepsilon
\eeno
for some $r>0$, then $(u,b)$ is regular at the point $z_0$.
\end{Proposition}

We also need the small energy interior regularity result in terms of the velocity only in \cite{WZ1}, and the according boundary regularity result see \cite{KK}.

\begin{Proposition}\label{prop:small regularity2}
Assume that $(u,b)$ is a suitable weak solution of (\ref{eq:MHD}) in $Q_1(z_0)$.
There exists an absolute constant $\varepsilon>0$ such that  if $u\in L^{p,q}$ near $z_0$ and
\ben\label{equ:2.6}
\limsup\limits_{r\rightarrow 0+} r^{-(\frac3p+\frac2q-1)}\|u\|_{L^{p,q}(Q_r(z_0))}<\varepsilon,
\een
with $p,q$ satisfying $1\leq \frac3p+\frac2q\leq 2$, $1\leq q\leq \infty$ and $(p,q)\neq (\infty,1)$.
Then $(u,b)$ is regular at the point $z_0$.
\end{Proposition}

Let $(u,b,\pi)$ be a solution of (\ref{eq:MHD}) and
introduce the following scaling: \ben\label{eq:scaling}
u^{\lambda}(x, t)={\lambda}u(\lambda x,\lambda^2 t),\quad
b^{\lambda}(x, t)={\lambda}b(\lambda x,\lambda^2 t),\quad
\pi^{\lambda}(x, t)={\lambda}^2\pi(\lambda x,\lambda^2 t), \een for any
$\lambda> 0,$
then the family $(u^{\lambda},b^{\lambda},\pi^{\lambda})$ is also a solution of
(\ref{eq:MHD}). For $z_0=(x_0,t_0)$, we define some invariant quantities under the scaling (\ref{eq:scaling}):

\beno
&&A(f,r)=\sup_{-r^2\leq t<0}r^{-1}\int_{B_r}|f(y,t)|^2dy,\quad\quad C(f,r)=r^{-2}\int_{Q_r}|f(y,s)|^3dyds,\\
&&E(f,r)=r^{-1}\int_{Q_r}|\nabla f(y,s)|^2dyds,\quad\quad K(f,r)=r^{-3}\int_{Q_r}|f(y,s)|^2dyds,
\eeno
for $f=u,b$ and
$$D(\pi,r)=r^{-2}\int_{Q_r}|\pi (y,s)|^{\f32}dyds,$$
$$\tilde{D}(\pi,r)=r^{-2}\int_{Q_r}|\pi (y,s)-(\pi)_{B_r}|^{\f32}dyds,\quad (\pi)_{B_r}=\frac{1}{|B_r|}\int_{B_r}\pi (y,s)dy.$$
Let $A(u,b;r)=A(u,r)+A(b,r)$, and $E(u,b;r)$, $C(u,b;r)$ and $K(u,b;r)$ denote similar notations.
We also introduce
\beno
&&G(f,p,q;r)=r^{1-\frac3p-\frac2q}\|f\|_{L^{p,q}(Q_r)},\quad \widetilde{G}(f,p,q;r)=r^{1-\frac3p-\frac2q}\|f-(f)_{B_r}\|_{L^{p,q}(Q_r)},\\
&&H(f,p,q;r)=r^{2-\frac3p-\frac2q}\|f\|_{L^{p,q}(Q_r)},\quad \widetilde{H}(f,p,q;r)=r^{2-\frac3p-\frac2q}\|f-(f)_{B_r}\|_{L^{p,q}(Q_r)}.
\eeno

Throughout this paper, we denote by $C_0$ a constant independent of $r,\rho$
and different from line to line.

\setcounter{equation}{0}
\section{Proof of Theorem \ref{thm:uh serrin}}

In this section, we'll prove Theorem \ref{thm:uh serrin} with the help of Theorem \ref{thm:uh interior}. Moreover, we assume that
\beno
\|u\|_{L^{\infty}_tL^3_x\big((-1,0)\times R^3\big)}+\|b(\cdot,-\frac12)\|_{L^3_x( R^3)}\leq C_0,
\eeno
which is reasonable from the assumptions of Theorem \ref{thm:uh serrin} and (\ref{eq:u b energy}).

First, we have the following embedding inequality and  Sobolev's interpolation inequality (for example, see \cite{CKN}):
\begin{Lemma}\label{lem:sobolev inequality}
i) For $2\leq \ell\leq 6$, $a=\frac34(\ell-2)$ and $f\in  H^1(R^3)$, we have
\ben\label{eq:embedding}
\int_{R^3}|f|^{\ell}\leq C_0\big(\int_{R^3}|\nabla f|^2\big)^a
\big(\int_{R^3}|f|^2\big)^{\frac \ell 2-a}.
\een
ii) For $f\in L^{\infty}(-1,0;L^2(R^3))\cap L^{2}(-1,0;\dot{H}^1(R^3)),$ we have
\ben\label{eq:embedding-parabolic}
\|f\|_{L^s_tL^l_x((-1,0)\times R^3)}\leq C_0 \|f\|_{L^{\infty}_tL^2_x((-1,0)\times R^3)}^{1-\frac2s} \|f\|_{L^2_t\dot{H}^1_x((-1,0)\times R^3)}^{\frac2s},
\een
where $\frac3l+\frac2s=\frac32$ with $2\leq s\leq \infty.$
\end{Lemma}

\begin{Lemma}\label{lem:b_3}
Under the assumption of  Theorem  \ref{thm:uh serrin}, we have
\beno
|b_3|^{\frac32}\in L^{\infty}_tL_x^{2}\cap L^2_t\dot{H}^1_x((-\frac12,0)\times\R^3).
\eeno
\end{Lemma}
\no{\bf Proof.}
Recall the third equation of the magnetic field:
\beno
\partial_tb_3-\Delta b_3+u\cdot \nabla b_3=b\cdot \nabla u_3 ,
\eeno
and multiplying $3|b_3|b_3$ on both sides of it, we have
\beno
\partial_t\big(|b_3|^3\big)-3\Delta b_3\big(|b_3|b_3\big)+u\cdot \nabla \big(|b_3|^3\big)=3b\cdot \nabla u_3 \big(|b_3|b_3\big).
\eeno
Integrating on $R^3$, using integration by parts and $\nabla\cdot u=0$ we derive that
\ben\label{eq:b3 1}
\partial_t\int_{R^3}|b_3|^3dx+\frac83\int_{R^3}|\nabla \big(|b_3|^{\frac32}\big)|^2dx=3\int_{R^3}\big(b\cdot \nabla u_3\big) \big(|b_3|b_3\big)dx\equiv 3I.
\een

Let $\nabla_h=(\partial_1,\partial_2)^T$ and $I=I_1+I_2,$ where
\beno
I_1=\int_{R^3}\big(b_h\cdot \nabla_h u_3\big) \big(|b_3|b_3\big)dx,\quad I_2=\int_{R^3}\big(b_3 \partial_3 u_3\big) \big(|b_3|b_3\big)dx.
\eeno
Obviously,
\beno
I_1&\leq& |\int_{R^3}\big(\nabla_h\cdot b_h\big)  u_3 \big(|b_3|b_3\big)dx|+2\int_{R^3}|b_h|  |u_3| |b_3||\nabla_h b_3|dx\\
&\leq& C_0\|u_3\|_{L^{3}(R^3)}\big[  \|\nabla_h\cdot b_h |b_3|^2\|_{L^{\frac32}(R^3)}+ \|b_3b_h|\nabla_h b_3| \|_{L^{\frac32}(R^3)}\big]\\
&\leq& C_0\big[  \||\nabla_h\cdot b_h| b_3^2\|_{L^{\frac32}(R^3)}+ \|b_3b_h|\nabla_h b_3| \|_{L^{\frac32}(R^3)}\big],
\eeno
and  the divergence-free property of $b$ implies that
\beno
I_2&\leq& C_0 \int_{R^3}|u_3| |\partial_3 b_3| |b_3|^2dx\\
&\leq& C_0\int_{R^3}|u_3| |\nabla_h\cdot b_h | b_3^2dx\leq C_0  \||\nabla_h\cdot b_h| b_3^2\|_{L^{\frac32}(R^3)}.
\eeno
From the above estimates of $I_1,I_2$ and (\ref{eq:b3 1}), we derive that
\ben\label{eq:b3 2}
\partial_t\int_{R^3}|b_3|^3dx+\frac83\int_{R^3}|\nabla \big(|b_3|^{\frac32}\big)|^2dx \leq C_0\big[  \||\nabla_h\cdot b_h| |b_3|^2\|_{L^{\frac32}(R^3)}+ \|b_3b_h|\nabla_h b_3| \|_{L^{\frac32}(R^3)}\big],
\een
and let
\beno
II_1\equiv \||\nabla_h\cdot b_h| b_3^2\|_{L^{\frac32}(R^3)},\quad II_2\equiv \|b_3b_h|\nabla_h b_3 |\|_{L^{\frac32}(R^3)}.
\eeno

{\bf Step I: Estimate of $II_1$.} Let $\frac3p+\frac2q=2$ with $1\leq q\leq \infty$ and
\beno
\|\nabla b_h\|_{L^{q}_tL^p_x((-\frac12,0)\times R^3)}\leq C_0.
\eeno
Then
by H\"{o}lder inequality we have
\beno
II_1&\leq& \|\nabla b_h\|_{L^p_x(R^3)}\|b_3\|_{L^{\frac{6p}{2p-3}}_x(R^3)}^2,
\eeno
where the end case $p=\infty$ or $p=\frac32$ still holds for the above inequality.  For any $\tau$ with $-\frac12<\tau<0$, integrating to time we have
\beno
\int_{-\frac12}^{\tau}II_1dt&\leq& \|\nabla b_h\|_{L^q_tL^p_x((-\frac12,\tau)\times R^3)}\|b_3\|^2_{L^{2q'}_tL^{\frac{6p}{2p-3}}_x((-\frac12,\tau)\times R^3)},
\eeno
where $\frac1q+\frac1{q'}=1.$

To apply Gronwall's inequality, we choose $q<\infty$ and by Lemma \ref{lem:sobolev inequality} $q',p$ should satisfy
\beno
\frac{2p-3}{2p}+\frac1{q'}=1,\quad 3\leq 2q'\leq \infty,
\eeno
which yields that $1\leq q\leq   3$ or $\frac94\leq p\leq \infty.$ Hence
\ben\label{eq:b3 3}
\|b_3\|_{L^{2q'}_tL^{\frac{6p}{2p-3}}_x((-\frac12,\tau)\times R^3)}\leq C_0\big[\|(|b_3|^{\frac32})\|_{L^{\infty}_{t}L^2_x((-\frac12,\tau)\times R^3)}+\|\nabla(|b_3|^{\frac32})\|_{L^2_{t,x}((-\frac12,\tau)\times R^3)}\big]^{\frac23}.
\een

{\bf Step II: Estimate of $II_2$.}

Let $\frac3l+\frac2s=1$ with $2\leq s\leq \infty$ and
\beno
\|b_h\|_{L^{s}_tL^l_x((-\frac12,0)\times R^3)}\leq C_0.
\eeno
Then
by H\"{o}lder inequality we have
\beno
II_2&\leq& C_0 \| b_h\|_{L^l_x(R^3)}\|\nabla(|b_3|^{\frac32})\|_{L^2(\R^3)}\|b_3\|_{L^{\frac{3p_1}{4}}_x(R^3)}^{\frac12},
\eeno
where $p_1$ satisfies
\beno
\frac3{2l}+\frac34+\frac1{p_1}=1,\quad 1\leq p_1\leq \infty.
\eeno
Integrating to time with the same $\tau$ as above, we have
\beno
\int_{-\frac12}^{\tau}II_2dt&\leq& \|b_h\|_{L^s_tL^l_x((-\frac12,\tau)\times R^3)} \|\nabla(|b_3|^{\frac32})\|_{L^2_{t,x}((-\frac12,\tau)\times R^3)}\|b_3\|_{L^{\frac{q_1}{2}}_tL^{\frac{3p_1}{4}}_x((-\frac12,\tau)\times R^3))}^{\frac12} ,
\eeno
where $q_1$ satisfies
\beno
\frac1{s}+\frac12+\frac1{q_1}=1,\quad 1\leq q_1\leq \infty.
\eeno

Using  $\frac3l+\frac2s=1$, we have
\beno
\frac4{p_1}+\frac4{q_1}=1,
\eeno
that is $|b_3|^{\frac32}\in L^{\frac{q_1}{3}}_tL^{\frac{p_1}{2}}_x$, and
\beno
3\frac2{p_1}+2\frac3{q_1}=\frac32,
\eeno
which yields that by Lemma \ref{lem:sobolev inequality}, for $2\leq \frac{q_1}{3}\leq \infty$,
\ben\label{eq:b3 4}
\|b_3\|_{L^{\frac{q_1}{2}}_tL^{\frac{3p_1}{4}}_x((-\frac12,\tau)\times R^3))}\leq C_0\big[\|(|b_3|^{\frac32})\|_{L^{\infty}_{t}L^2_x((-\frac12,\tau)\times R^3)}+\|\nabla(|b_3|^{\frac32})\|_{L^2_{t,x}((-\frac12,\tau)\times R^3)}\big]^{\frac23}.
\een
Then $s$ must be $2\leq s\leq 3$ or $9\leq l\leq \infty.$

{\bf Step III: Arguments.} From the above two estimates, we obtain that if the condition $i)$ of (\ref{eq:bh condition}) holds, i.e.
\beno
\|\nabla b_h\|_{L^{q}_tL^p_x((-\frac12,0)\times R^3)}\leq C_0, \quad \frac94\leq p<3,
\eeno
then the embedding inequality implies
\beno
\|b_h\|_{L^{q}_tL^l_x((-\frac12,0)\times R^3)}\leq C_0, \quad \frac3l+\frac2q=1,\quad 9\leq l<\infty,
\eeno
thus the estimates (\ref{eq:b3 3})-(\ref{eq:b3 4}) hold, which yields that for $-\frac12<\tau<0$, there holds
\ben\label{eq:b3 2}
&&\int_{-\frac12}^{\tau}\int_{R^3}\partial_t(|b_3|^3)dxdt+\frac83\int_{-\frac12}^{\tau}\int_{R^3}|\nabla \big(|b_3|^{\frac32}\big)|^2dxdt \nonumber\\
&\leq& C_0\big[\|(|b_3|^{\frac32})\|_{L^{\infty}_{t}L^2_x((-\frac12,\tau)\times R^3)}+\|\nabla(|b_3|^{\frac32})\|_{L^2_{t,x}((-\frac12,\tau)\times R^3)}\big]^{\frac43}.
\een
Note that $\|b(\cdot,-\frac12)\|_{L^3_x( R^3)}\leq C_0$. Hence, Gronwall's inequality implies
\beno
\|(|b_3|^{\frac32})\|_{L^{\infty}_{t}L^2_x((-\frac12,\tau)\times R^3)}+\|\nabla(|b_3|^{\frac32})\|_{L^2_{t,x}((-\frac12,\tau)\times R^3)}<\infty,
\eeno
and taking the supremum of $\tau\rightarrow 0$, we have
\ben\label{eq:b3 estimate}
\|(|b_3|^{\frac32})\|_{L^{\infty}_{t}L^2_x((-\frac12,0)\times R^3)}+\|\nabla(|b_3|^{\frac32})\|_{L^2_{t,x}((-\frac12,0)\times R^3)}\leq C_0.
\een

On the other hand, when the condition $ii)$ of (\ref{eq:bh condition}) holds, simialr estimates hold. The lemma is proved.\endproof

{\bf Proof of Theorem \ref{thm:uh serrin}:} Due to the assumptions on $b_h$, we have
\beno
\|b_h\|_{L^{s}_tL^l_x((-\frac12,0)\times R^3)}\leq C_0, \quad \frac3l+\frac2s=1,\quad 9\leq l\leq \infty.
\eeno
Moreover, Lemma 3.2 yields that
\beno
\|b_3\|_{L^{5}_{t,x}((-\frac12,0)\times R^3)}<\infty.
\eeno
Obviously, the conditions of Theorem \ref{thm:uh interior} are satisfied, thus the proof is complete.\endproof

\setcounter{equation}{0}
\section{Blow-up analysis and Proof of Theorem \ref{thm:uh interior}}
We will apply Proposition \ref{prop:small regularity} to prove the interior regularity of the solution by blow-up analysis, which was early used for the 3D Navier-Stokes equations in \cite{Lin}, see also \cite{Se2}. Note that the velocity $u$ is in the critical class, hence backward uniqueness results in \cite{ESS} are still needed. Firstly, we prove that the basic energy norms $A(u,b;r)$, $E(u,b;r)$, and $\tilde{D}(\pi,r)$ are uniformly bounded for all $0<r<1.$ (see Theorem \ref{thm:bounded A E}); secondly, a standard compactness argument for suitable weak solutions of the 3-D MHD equations and backward uniqueness results imply Theorem \ref{thm:uh interior}.

\subsection{Bounded estimates of $A(u,b;r)$ and $E(u,b;r)$}

To ensure the validness of blow-up analysis, we have to prove that $A(u,b;r)$, $E(u,b;r)$, and $\tilde{D}(\pi,r)$ are uniformly bounded for all $0<r<r_1$ with some $r_1>0.$
\begin{Theorem}\label{thm:bounded A E}
Under the assumptions of \ref{thm:uh interior}, there exists a $r_1>0$ such that
\ben\label{eq:u_h nablau_h}
A(u,b;r)+E(u,b;r)+\tilde{D}(\pi,r)<\infty,\quad \quad 0<r<r_1,
\een
where $r_1$ depends on $C(u,b;1)$ and $\tilde{D}(\pi,1)$.
\end{Theorem}

For completeness, we supply the following technical lemmas. First, we will control $A(u,b;r)+E(u,b;r)$ in terms of the other scaling invariant quantities by using the following local inequality.
\begin{Lemma}\label{lem:local energy}
Let $0<4r<\rho<r_0$ and $1\le p, q\le \infty$. There holds
\beno
&&A(u,b;r)+E(u,b;r)\nonumber\\
&&\leq
C_0\big(\frac{r}{\rho}\big)^2K(u,b;\rho)+C_0\big(\frac{\rho}{r}\big)^2
\Big[C(u,\rho)+C(u,\rho)^{1/3}\big(C(b,\rho)^{2/3}+\tilde{D}(\pi,\rho)^{2/3}\big)\Big].
 \eeno
\end{Lemma}

\no{\bf Proof.}\,Let $\zeta$ be a cutoff function, which vanishes
outside of $Q_{\rho}$ and equals 1 in $Q_{\rho/2}$, and satisfies
$$|\nabla\zeta|\leq C_0\rho^{-1},\quad |\partial_t\zeta|, |\triangle\zeta|\leq C_0\rho^{-2}.$$
Define the backward heat kernel as
$$\Gamma(x,t)=\frac{1}{4\pi(r^2-t)^{\frac32}}e^{-\frac{|x|^2}{4(r^2-t)}}.$$
Taking the test function $\phi=\Gamma\zeta$ in the local energy inequality, and noting $(\partial_t+\triangle)\Gamma=0$,
we have
\beno
&&\sup_t\int_{B_{\rho}}(|u|^2+|b|^2)\phi dx+\int_{Q_{\rho}}(|\nabla u|^2+|\nabla b|^2)\phi dxdt\\
&&\leq \int_{Q_{\rho}}\big[(|u|^2+|b|^2)(\Gamma\triangle\zeta+\Gamma\partial_t\zeta+2\nabla\Gamma\cdot\nabla\zeta)
+|\nabla\phi||u|(|u|^2+|b|^2+2|\pi-\pi_{B_{\rho}}|)\big]dxdt.
\eeno
Some direct computations imply
\beno
&&\Gamma(x,t)\geq C_0^{-1}r^{-3}\quad {\rm in} \,\, Q_r;\\
&&|\nabla\phi|\leq |\nabla\Gamma|\zeta+\Gamma|\nabla\zeta|\leq C_0r^{-4};\\
&&|\Gamma\triangle\zeta|+|\Gamma\partial_t\zeta|+2|\nabla\Gamma\cdot\nabla\zeta|\leq C_0\rho^{-5},
\eeno
from which and H\"{o}lder inequality, Lemma \ref{lem:local energy} follows.\endproof

The following is an interpolation inequality.

\begin{Lemma}\label{lem:C(u;r)}
For any $0<r<r_0$, let $\frac3p+\frac2q=2$ with $1\leq q\leq \infty$, there holds
\beno
C(f,r)\leq C_0G(f,p,q;r)\big(E(f,r)+A(f,r)\big),
\eeno
where $f=u,b$.
\end{Lemma}

\no{\bf Proof.} Without loss of generality, we consider the estimate of $u$.
By H\"{o}lder inequality and Sobolev inequality, we get
\beno
\int_{B_r}|u|^3dx &= &\int_{B_r}|u|^{3\alpha+3\beta+3-3\alpha-3\beta}dx\\
&\leq & \big(\int_{B_r}|u|^2dx\big)^{3\alpha/2}\big(\int_{B_r}|u|^6dx\big)^{\beta/2}
\big(\int_{B_r}|u|^pdx\big)^{{(3-3\alpha-3\beta)}/p}\\
&\leq & C_0\big(\int_{B_r}|u|^2dx\big)^{3\alpha/2}\big(\int_{B_r}|\nabla u|^2+|u|^2dx\big)^{3\beta/2}
\big(\int_{B_r}|u|^pdx\big)^{{(3-3\alpha-3\beta)}/p},
\eeno
where  $\alpha, \beta$ are chosen so that
\beno
\frac13=\frac{\alpha}{2}+\frac{\beta}{6}+\frac{1-\alpha-\beta}{p},\quad
\frac{3\beta}{2}+\frac{3-3\alpha-3\beta}{q}=1.
\eeno
Taking $\alpha=\frac{2p-3}{3p}$ and $\beta=\frac1p$, we get
 \beno
\int_{Q_r}|u|^3dx
&\leq& C_0\big(\sup_{-r^2<t<0}\int_{B_r}|u|^2dx\big)^{1-\frac{3}{2p}}
\big(\int_{Q_r}|\nabla u|^2+|u|^2dxdt\big)^{\f 3{2p}}\\
&&\quad\times \Big(\int_{-r^2}^0\big(\int_{B_r}|u|^pdx\big)^{\f q p}dt\Big)^{\f 1 q},
\eeno
which yields the required inequality.\endproof

We present the estimate of the pressure in terms of scaling invariant quantities, see also \cite{Se1}.

\begin{Lemma}\label{lem:pressure}
Let $(u,b)$ be a suitable weak solution of (\ref{eq:MHD}) in $Q_1$. Then there hold
\ben\label{eq:pressure}
\widetilde{D}(\pi,r)\leq C_0\big((\frac{r}{\rho})^{5/2}\widetilde{D}(\pi,\rho)+(\frac{\rho}{r})^{2}C(u,b;\rho)\big),
\een
for any $0<4r<\rho<1.$
\end{Lemma}

\no{\bf Proof.}
Note that $\pi$ satisfies the following equation in distribution sense:
$$-\triangle\pi=\partial_i\partial_j(\hat{u}_i\hat{u}_j-\hat{b}_i\hat{b}_j),$$
where $\hat{u}=u-(u)_{B_{\rho}}$ and $\hat{b}=b-(b)_{B_{\rho}}$. Let $\zeta$ be a cut-off function,
which equals 1 in $Q_{\rho/2}$ and vanishes outside of $Q_{\rho}$.
Set $\pi=\pi_1+\pi_2$ with
$$
\pi_1=\frac{1}{4\pi}\int_{\R^3}\frac{1}{|x-y|}
[\partial_i\partial_j\big((\hat{u}_i\hat{u}_j-\hat{b}_i\hat{b}_j)\zeta^2\big)],
$$
and $\pi_2$ is harmonic in $Q_{\rho/2}$.

Due to the Calderon-Zygmund inequality, we have
$$
\int_{B_{\rho}}|\pi_1|^{\frac32}dx\leq C_0\int_{B_{\rho}}|\hat{u}|^{3}+|\hat{b}|^{3}dx.$$
Since $\pi_2$ is harmonic in $Q_{\rho/2}$, we have
\beno
\int_{B_{r}}|\pi_2-(\pi_2)_{B_r}|^{\frac32}dx&\leq& C_0r^{3+\frac32}\sup_{B_{\rho/4}}|\nabla\pi_2|^{\frac32}\\
&\leq& C_0\big(\frac{r}{\rho}\big)^{3+\frac32}\int_{B_{\rho/2}}|\pi_2-(\pi_2)_{B_{\rho/2}}|^{\frac32}dx.
\eeno
Hence we infer that
\beno
\tilde{D}(\pi,r)
\leq \tilde{D}(\pi_1,r)+\tilde{D}(\pi_2,r)
\leq C_0\big((\frac{r}{\rho})^{5/2}\widetilde{D}(\pi,\rho)+(\frac{\rho}{r})^{2}C(u,b;\rho)\big).
\eeno
\endproof\vspace{0.1cm}

{\bf Proof of  Theorem \ref{thm:bounded A E}}. Without loss of generality, we assume that
$$\sup_{0<r<1}G(b_h,p,q;r)+\sup_{0<r<1}G(b_3,l,s;r)\leq C_0,\quad $$
for some $(p,q)$ and $(l,s)$ satisfying
$$\frac3p+\frac2q=2,\quad \frac3l+\frac2s=2.$$
 Then, by Lemma \ref{lem:C(u;r)} we have
$$C(b,\rho)\leq C_0[C(b_h,\rho)+C(b_3,\rho)]\leq C_0[A(b,\rho)+E(b,\rho)]. $$
Thus, by the local energy inequality and $C(u,\rho)\leq C_0$, we have
\beno
&&A(u,b;r)+E(u,b;r)\nonumber\\
&&\leq
C_0\big(\frac{r}{\rho}\big)^2A(u,b;\rho)+C_0\big(\frac{\rho}{r}\big)^2
\Big[1+\big(C(b,\rho)^{2/3}+\tilde{D}(\pi,\rho)^{2/3}\big)\Big]\\
&&\leq
C_0\big(\frac{r}{\rho}\big)^2[A(u,b;\rho)+E(u,b;\rho)]+C_0\big(\frac{\rho}{r}\big)^2
\Big[\big(\frac{\rho}{r}\big)^8+\tilde{D}(\pi,\rho)^{2/3}\big)\Big].
 \eeno
 Take $0<8r<\rho<r_0$ and set
\beno
F(r)=A(u,b;r)+E(u,b;r)+\varepsilon^{-1/2} \tilde{D}(\pi,r)^{2/3}.
\eeno
Due to Lemma \ref{lem:pressure}, we get
\beno
\tilde{D}(\pi,r)&\leq& C_0\big((\frac{r}{\rho})^{5/2}\widetilde{D}(\pi,\rho)+(\frac{\rho}{r})^{2}C(u,b;\rho)\big)\\
&\leq& C_0(\frac{r}{\rho})^{5/2}\widetilde{D}(\pi,\rho)+C_0(\frac{\rho}{r})^{2}\big(1+A(b,\rho)+E(b,\rho)\big)
\eeno
Hence,
\beno
F(r)\leq C_0\big[\big(\frac{r}{\rho}\big)^2+ \big(\frac{\rho}{r}\big)^2\varepsilon^{1/2}+(\frac{r}{\rho})^{5/3} \big]F(\rho)+C_0\varepsilon^{-3/2}\big(\frac{\rho}{r}\big)^{10},
\eeno
and choosing $\theta=\big(\frac{\rho}{r}\big)$ and $\varepsilon$ sufficiently small, we get
\beno
F(\theta\rho)\leq \frac12F(\rho)+C_0\varepsilon^{-3/2}\theta^{-10},
\eeno
which yields that there exists a $r_1>0$ such that
\beno
\sup_{0<r<r_1}F(r)<C_1,
\eeno
where $C_1$ depends on $C_0$, $C(u,b;1), \tilde{D}(\pi,1).$  \endproof\vspace{0.1cm}

\subsection{Proof of Theorem \ref{thm:uh interior} }

The proof of Theorem \ref{thm:uh interior} is based on the blow-up analysis and unique continuation theorem, for example see \cite{Se2}, \cite{ESS}.

We assume that $\|u\|_{L^{\infty}_tL^3_x(Q_1)}\leq C_0$, and
\beno
&&i)\liminf_{r\rightarrow 0}G(b_h,p,q;r)=0,\quad \limsup_{0<r<1}G(b_h,p,q;r)\leq C_0, \quad \frac3p+\frac2q=2,\quad 1\leq p\leq \infty;\\
&&ii)\sup_{0<r<1}G(b_3,l,s;r)<C_0,\quad \frac3l+\frac2s=2,\quad 1\leq s\leq \infty.
\eeno
Then, by the local energy inequality in Proposition \ref{prop:local} and Theorem \ref{thm:bounded A E}, we have
\beno
A(u,b;r)+E(u,b;r)+\tilde{D}(\pi,r)\leq C_1, \quad {\rm for \,\, all} \,\,0<r<1,
\eeno
where $C_1>0$ may depend on $C(u,b;1)$ and $D(\pi,1)$. Moreover, suppose that
$C(u,b;1)+D(\pi,1)\leq C_0$, which is reasonable by the definition of suitable weak solutions.

Suppose that the statement of the theorem is false. Then there exist
a series of suitable weak solutions $(v^k,\bar{b}^k,\bar{\pi}^k)$ and
 $r_k\downarrow 0$ such that
\ben\label{eq:A Ek }
A(v^k,\bar{b}^k;r)+E(v^k,\bar{b}^k;r)+\tilde{D}(\bar{\pi}^k,r)\leq C_1, \quad {\rm for \,\, all}\,\, 0<r<1,
\een
and
\beno
G(\bar{b}_h^k,p,q;r_k)\rightarrow0,\quad {\rm as}\,\,r_k\rightarrow 0.
\eeno
Moreover, $(0,0)$ is a singular point of $(v^k,\bar{b}^k,\bar{\pi}^k)$.

We denote
\beno
u^k(y,s)=r_kv^k(r_ky,r_k^2s),\quad b^k(y,s)=r_k\bar{b}^k(r_ky,r_k^2s),\quad \pi^k(y,s)=r_k^2 \bar{\pi}^k(r_ky,r_k^2s),
\eeno
where $(y,s)\in B_{\frac1{r_k}}\times (-\f1{r_k^2},0).$
Then it follows from (\ref{eq:A Ek }) that
\ben\label{eq:A Ek2}
&&A(u^k,b^k;r)+E(u^k,b^k;r)+\tilde{D}(\pi^k,r)\leq C_1, \quad {\rm for \,\, all}\,\, 0<r<1,\nonumber\\
&& G(b_h^k,p,q;1)\rightarrow0, \quad {\rm as}\,\,r_k\rightarrow 0.\\
&& \|u^k\|_{L^{\infty}_tL^3_x(B_{\frac1{r_k}}\times (-\f1{r_k^2},0))}\leq C_0.\nonumber
\een
For any $a,T>0$, choose sufficiently large $k$ such that
\ben\label{eq:4.5}
&&\|u^k\|_{L^{\infty,2}((-T,0)\times B_a)}+\|b^k\|_{L^{\infty,2}((-T,0)\times B_a))}\nonumber\\
&&\quad+\|\nabla u^k\|_{L^{2}((-T,0)\times B_a))}
+\|\nabla b^k\|_{L^{2}((-T,0)\times B_a))}\leq c(a,T).
\een
Hence, $u^k\cdot\nabla u^k,u^k\cdot\nabla b^k,b^k\cdot\nabla u^k, b^k\cdot\nabla b^k\in L_t^{\frac32}L_x^{\frac98}(Q_a)$. This gives by the linear Stokes theory \cite{ESS} that
\beno
|\partial_t u^k|+|\Delta u^k|+|\partial_t b^k|+|\Delta b^k|+|\nabla p^k|\in L_t^{\frac32}L_x^{\frac98}(Q_{3a/4}).
\eeno

Then Lions-Aubin's lemma ensures that there exists $(u,b,\pi)$ such that for any $a,T>0$ (up to subsequence),
\beno
&&u^k\rightarrow u ,\,b^k\rightarrow b, \quad {\rm in} \quad L^3((-T,0)\times B_a),\\
&&u^k\rightarrow u, \,b^k\rightarrow b, \quad {\rm in} \quad C([-T,0];L^{9/8}(B_{a})),\\
&&\pi^k\rightharpoonup \pi\quad {\rm in} \quad L^\f {3}2((-T,0)\times B_a),\\
&&\|u\|_{L^{\infty}_tL^3_x((-T,0)\times R^3)}\leq C_0,\\
&& b_h^k\rightharpoonup b_h=0, \quad {\rm in} \quad L^q((-T,0);L^p( B_a)),
\eeno
as $k\rightarrow+\infty$.

Hence $\partial_3b_3=0$ and $b\cdot\nabla b=0$ due to the velocity field equations, and we get
\begin{equation}\label{eq:u}
u_t-\Delta u+u\cdot \nabla u=-\nabla \pi,\quad \nabla\cdot u=0.
\end{equation}
Using the property of weak convergence, by (\ref{eq:A Ek2}) we have
$$C(u,1)+\tilde{D}(\pi,1)\leq C_0,$$
and due to
$u\in L^{\infty,3}_{t,x}$ , the well-known result in \cite{ESS} yields that $\|u\|_{L^{\infty}(Q_{\frac12})}\leq C_0.$ Applying the interior regularity criteria in
Proposition \ref{prop:small regularity2}, we obtained that $\|b\|_{L^{\infty}(Q_{\frac12})}\leq C_0.$

Since $(0,0)$ is a singular point of $(v^k,\bar{b}^k,\bar{\pi}^k)$, by small regularity results in Proposition \ref{prop:small regularity} we have
\beno
\varepsilon<C(v^k,r)+C(\bar{b}^k,r)+\tilde{D}(\bar{\pi}^k,r),
\eeno
for any $0<r<1.$ Thus,
\beno
\varepsilon<C(u^k,r)+C(b^k,r)+\tilde{D}(\pi^k,r),
\eeno
for any $0<r<1.$

Take the supremum limit as $k\rightarrow\infty$, we have
\beno
\varepsilon<C_0r^3+\tilde{D}(\pi^k,r),
\eeno
for any $0<r<1.$

By the pressure estimate in Lemma \ref{lem:pressure} and (\ref{eq:A Ek2}), we have
\beno
\widetilde{D}(\pi^k,r)\leq C_0\big((\frac{r}{\rho})^{5/2}\widetilde{D}(\pi^k,\rho)+(\frac{\rho}{r})^{2}C(u^k,b^k;\rho)\big)\leq C_0(\frac{r}{\rho})^{5/2}+C_0(\frac{\rho}{r})^{2}C(u^k,b^k;\rho),
\eeno
for any $0<r<\rho<1.$ Choose $\rho=\sqrt{r}$, then
\beno
\limsup_{k\rightarrow\infty}\widetilde{D}(\pi^k,r)\leq C_0 \sqrt{r},
\eeno
for any $0<r<1.$

Hence, we have
$
\varepsilon<C_0r^3+C_0 \sqrt{r},
$
for any $0<r<1.$ Obviously, it's a contradiction. The proof is complete.
\endproof

\bigskip

\noindent {\bf Acknowledgments.}
The author would like to thank Prof. Zhifei Zhang for some helpful discussions. W. Wang was supported by NSFC 11301048, "the Fundamental Research Funds for the Central Universities" and The Institute of Mathematical Sciences of CUHK.




\begin{thebibliography}{WWW}

\bibitem{CKN} L. Caffarelli, R. Kohn and L. Nirenberg,
{\it Partial regularity of suitable weak solutions of the Navier-Stokes equations},
Comm. Pure Appl. Math., 35(1982), 771-831.

\bibitem{Cao} C. Cao and J. Wu, {\it Two regularity criteria for the 3D MHD equations},
J. Differential Equations, 248(2010), 2263-2274.

\bibitem{CMZ1} Q. Chen, C. Miao and Z. Zhang,
{\it The Beale-Kato-Majda  criterion to the 3D Magneto-hydrodynamics equations},
Comm. Math. Phys., 275(2007), 861-872.

\bibitem{CMZ2} Q. Chen, C. Miao and Z. Zhang,
{\it On the regularity criterion of weak solution for the 3D viscous Magneto-hydrodynamics equations},
Comm. Math. Phys., 284(2008), 919-930.

\bibitem{DD} H. Dong and D. Du, {\it The Navier-Stokes equations in the critical Lebesgue space,}
Comm. Math. Phys., 292 (2009), no. 3, 811-827.

\bibitem{DL} G. Duvaut and J.-L. Lions,
{\it In\'{e}quations en thermo\'{e}lasticit\'{e} et magn\'{e}tohydrodynamique,}
Arch. Rational. Mech. Anal., 46(1972), 241-279.

\bibitem{ESS} L. Escauriaza, G. A. Seregin and V. \v{S}ver\'{a}k,
{\it $L^{3,\infty}$ solutions to the Navier-Stokes equations and backward uniqueness},
Russian Math. Surveys, 58(2003), 211-250.

\bibitem{HX1} C. He and Z. Xin, {\it On the regularity of weak solutions to the magnetohydrodynamic equations,}
J. Differential Equations, 213(2005), 235-254.

\bibitem{HX2} C. He and Z. Xin, {\it Partial regularity of suitable weak solutions to the incompressible magnetohydrodynamic equations,} J. Funct. Anal., 227(2005), 113-152.

\bibitem{KL} K. Kang and J. Lee, {\it Interior regularity criteria for suitable weak solutions of the magnetohydrodynamics equations},
J. Differentional Equations, 247(2009), 2310-2330.

%



\bibitem{MNS} A. Mahalov, B. Nicolaenko and T. Shilkin, {\it $L^{3,\infty}$-solutions to the MHD equations,}
Zap. Nau\v{c}n. Sem. S.-Peterburg. Otdel. Mat. Inst. Steklov. (POMI) 336(2006), 112-132.


\bibitem{PP} H. Politano, A. Pouquet  and P.-L. Sulem, {\it Current and vorticity dynamics in three-dimensional magnetohydrodynamic turbulence}, Phys. Plasmas, 2(1995), 2931-2939.


\bibitem{Se1} G. A. Seregin, {\it Estimate of suitable solutions to the Navier-Stokes equations in critical Morrey spaces}, Journal of Mathematical Sciences, 143(2007),  2961-2968.

\bibitem{Se2} G. A. Seregin, {\it On the local regularity of suitable weak solutions of the Navier-Stokes equations.} (Russian) Uspekhi Mat. Nauk 62 (2007), no. 3(375), 149-168; translation in Russian Math. Surveys 62 (2007), no. 3, 595-614.

\bibitem{ST} M. Sermange and R. Teman, {\it Some mathematical questions related to the MHD equations,}
Comm. Pure Appl. Math., 36(1983), 635-664.



\bibitem{WZ1} W. Wang and Z. Zhang, {\it On the interior regularity criteria for suitable weak solutions of the Magneto-hydrodynamics equations}, SIAM J. Math. Anal., 45 (2013), no. 5, 2666-677.

\bibitem{WZ2} W. Wang and Z. Zhang, {\it Limiting case for the regularity criterion to the 3-D magneto-hydrodynamics equations}, J. Differential Equations, 252 (2012), no. 10, 5751-5762.

\bibitem{Wu1} J. Wu, {\it Bounds and new approaches for the 3D MHD equations},
J. Nonlinear Sci., 12(2002), 395-413.

\bibitem{Wu2} J. Wu, {\it Regularity results for weak solutions of the 3D MHD equations},
Discrete Contin. Dyn. Syst., 10(2004), 543-556.



\bibitem{Zhou1} Y. Zhou, {\it Remarks on regularities for the 3D MHD equations}, Discrete. Contin. Dyn. Syst., 12(2005), 881-886.


\end{thebibliography}
\end{document}